\documentclass[12pt]{article}        

\usepackage{color,graphicx}
\usepackage{latexsym,amssymb,amsmath}

\newtheorem{theorem}{Theorem}
\newtheorem{remark}{{\it{Remark}}}

\begin{document}


\title{A general formulation of Non-Local Dirichlet forms on infinite dimensional topological 
vector spaces and its applications, and corresponding subjects}
\author{
{{Sergio
 \textsc{Albeverio}}}
          \footnote{Inst. Angewandte Mathematik, 
and HCM, Univ. Bonn, Germany, 
email \, :albeverio@iam.uni-bonn.de},
\quad 
{{Toshinao \textsc{Kagawa}}} \footnote{Dept. Information Systems Kanagawa Univ., Yokohama, Japan }, 
\quad
{{Shyuji \textsc{Kawaski}}} \footnote{Dept. Mathematical Sciences and Phyisics,  Iwate Univ., Morioka, Japan }, 
\\
 {{Yumi \textsc{Yahagi}}} \footnote{Dept. Mathematical information Tokyo Univ. of Information, Chiba, Japan}, 
\quad
 and \quad
{{Minoru \textsc{W. Yoshida}}} \footnote{Dept. Information Systems Kanagawa Univ., Yokohama, Japan,  
 email:\, washizuminoru@hotmail.com }
 }

\date{{\textcolor{blue}{SEMINAR GFM-GRUPO DE F{I'}SICA MATEM{\'A}TICA}}\\ Universidada de Lisboa\, 13rd Sept. 2023} 

\maketitle

\begin{abstract}
Though the present talk we give a concise explanation on  
the results derived and considered in the following two papers. 
Namely, here, we  break the following subjects  down:

\medskip

\textcolor{blue}{
i)  Open Access
Published: 14 August 2021
Non-local Markovian Symmetric Forms on Infinite Dimensional Spaces I.
The closability and quasi-regularity.
By Sergio Albeverio, Toshinao Kagawa, Yumi Yahagi,  Minoru W. Yoshida,
Communications in Mathematical Physics volume 388, pages659–706 (2021).
}

\textcolor{blue}{
ii)  Open Access
Published: 16 August 2022
Non-local Markovian Symmetric Forms on Infinite Dimensional Spaces
Part 2. Applications: Non Local Stochastic Quantization of Space Cut-Off Quantum Fields and Infinite Particle Systems.
By Sergio Albeverio, Toshinao Kagawa, Shuji Kawasaki, Yumi Yahagi,  Minoru W. Yoshida,
Potential Analysis (2022).
}
\bigskip

\medskip

A main purpose of this talk is an explanation on the abstract framework of 
a formulation of non-local Markovian forms defined on the probability spaces 
on (separable) infinite dimensional topological vector spaces, and the stochastic quantizations.

A little precisely, 
general theorems on the closability and quasi-regularity of non-local Markovian symmetric forms on probability  spaces $(S, {\cal B}(S), \mu)$, with $S$ Fr{\'e}chet spaces such that $S \subset  
{\mathbb R}^{\mathbb N}$, ${\cal B}(S)$ is the Borel $\sigma$-field of $S$, and $\mu$ is a Borel  probability measure on $S$, are introduced. Precisely, a family of non-local Markovian symmetric forms  
${\cal E}_{(\alpha)}$,$0 < \alpha < 2$, acting in each given $L^2(S; \mu)$ is defined, the index  
$\alpha$ characterizing the order of the non-locality.

Then, we see that all the forms ${\cal E}_{(\alpha)}$ defined on $\bigcup_{n \in {\mathbb N}}  C^{\infty}_0({\mathbb R}^n)$ are closable in $L^2(S;\mu)$.
Moreover, sufficient conditions under which the closure of the closable forms, that are Dirichlet forms,  
become strictly quasi-regular, are given.
Through these, an existence theorem for Hunt processes properly associated to the Dirichlet forms
is given.

The application of the above general theorems we consider the problem of stochastic quantizations 
of Euclidean $\Phi^4_d$ fields, for $d =2, 3$, by means of these Hunt processes is indicated. 
Also, some of the important corresponding problems, e.g., the stochastic mechanics, Malliavin calculus for jump processes, etc., shall be mentioned.


\end{abstract}

\tableofcontents

\section{What is the stochastic quantization}

{\textcolor{blue}{  \subsection{An elementary concept corresponding to Markov chains}}}
\medskip
\medskip

We recall the corresponding 
problem defined in the framework of Markov chains with finite state spaces.
Let $n \in {\mathbb N}$, suppose that we are given an $n\times n$ Markovian matrix $M \equiv (p_{ij})$, $i,j=1, \dots, n$. Thus the state space $S$ of the Markov process in consideration is $S = \{1, \dots, n\}$. 
Then a probability distribution $\mu = (p_1, \dots, p_n)$, $0 \leq p_i \leq 1$, 
$\sum_{i=1}^n p_i = 1$, that satisfies $\mu \, M = \mu$ is an invariant measure 
of the Markov process defined through $M$, i.e., $\mu$ is an {\it{eigenvector}} of 
$M$ with the {\it{eigenvalue}}  1.

 Conversely, suppose that we are given a 
 probability distribution $\mu = (p_1, \dots, p_n)$, $0 \leq p_i \leq 1$, 
$\sum_{i=1}^n p_i = 1$, and suppose that we are asked to find an $n\times n$ Markovian matrix $M$ such that 
$$
{\textcolor{blue}{
\mu \, M = \mu}}. 
\eqno{(*)} $$
The latter problem is the
 {\textcolor{blue}{{\it{stochastic quantization}} }}
of the probability measure $\mu$. It is the interpretation of 
the  {\it{stochastic quantization}} in the framework of the Markov process with finite state 
space. Obvously, in general there  exist many $M$ that satisfy $\mu \, M = \mu$ for a 
given $\mu$.  
 {\it{Ergodicity}} of  the associated Markov process is related to uniqueness of $M$.

{\textcolor{blue}{\subsection{ Stochastic quantization by means of the Markovian symmetric forms}}}
Under the situation of the previous subsection with the state space $S=\{1, \cdots, n\} \subset {\mathbb R}^1$, if we find the $M$ satisfying ($\ast$), then for any real $n$ vector 
$\kappa = (k_1, \dots, q_n)$, 
$$(\mu  \, M^m)  \cdot \kappa  = \mu  \cdot \kappa, \qquad \forall  m \in {\mathbb N}.$$ 
In particular, if $M$ is a symmetric (Markovian) matrix, then,
$${\textcolor{blue}{
(\mu  \, M^m)  \cdot \kappa  = \mu  \cdot \kappa 
= \mu \cdot (\kappa M^m), 
\qquad \forall  m \in {\mathbb N}.
}}
\eqno{(**)} $$

As the analogy, let $S = {\mathbb R}^1$ and $\mu$ is a probability measure on $S$, if we have 
a symmetric Markovian semi-graup $e^{-Ht}$, $t \geq 0$ (a particular case of $L^2_{\mu}$ contraction semi-group), then

the corresponding formula to ($**$) becomes, by the symmetry, 
$${\textcolor{blue}{
\int_{\mathbb R} (e^{-Ht} f)(x) \, \mu(dx)
= \int_{\mathbb R} f(x) \, (e^{-Ht} 1)(x) \, \mu(dx)}}
$$
$${\textcolor{blue}{ = \int_{\mathbb R} f(x) \, \mu(dx), }}
$$
$${\textcolor{blue}{ 
{\mbox{for any real-valued bounded measurable $f$ and \, \, $\forall t \geq 0$}}.
}}
\eqno{(***)}$$
By the differentiation, ($***$) holds, if $e^{-Ht} H 1 =0$ or equivalently $H 1 =0$, i.e., 
the function \, $1$\,  is the eigenvector of $H$ with the eigenvalue $0$.

Hence, suppose that a symmetric form ${\cal E} (\cdot, \cdot)$ is the one corresponding to $H$, 
$${\cal E} (f, g) = \int_{\mathbb R} (Hf)(x) \, g(x) \, \mu(dx),$$
then the stochastic quantization of $\mu$, i.e., to identify $e^{- Ht}$, equivalently, to identify 
$H$ or ${\cal E}(\cdot, \cdot)$ is 
interpreted as the problem as follows: \\
Define a symmetric Markovian form (the Dirichlet form) 
${\cal E} (\cdot, \cdot)$ on $L^2_{\mu}$ such that 
$$
{\textcolor{blue}{
{\cal E}(f, 1) = 0, \qquad   \forall f \in {\cal D}({\cal E}).
}} \eqno{(****)}
$$

{\textcolor{blue}{\subsection{Classical (local, i.e. diffusion) Dirichlet forms\\
 on real  $f(x)$, $x \in 
S \subset {\mathbb R}^m$, $m < \infty$,\\
 i.e. {\it{finite}} dimensional topological vector spaces}}}
For simplicity, let $m = 1$, and let $\rho \in C^{\infty} ({\mathbb R} \to {\mathbb R}_+)$, which is a strictly positive probability density function. 
For $f,\, g \in C^{\infty}_0( {\mathbb R} \to {\mathbb R})$ define a symmetric form 
${\cal E}(\cdot, \cdot)$ 
on $L^2_{\rho dx}$
$${\cal E}(f, g) = \int_{\mathbb R} \frac{d}{dx}f(x) \cdot  \frac{d}{dx}g(x) \cdot \rho(x) dx$$
$$ = \int_{\mathbb R} f(x)  \cdot \{ (- \frac{d^2}{dx^2} - \frac{\rho'}{\rho} \frac{d}{dx})g(x) \}  \cdot 
\rho (x) \, dx.$$

Then (see ($****$))  a closed extension of ${\cal E}$ defines a symmetric Markovian semi-group to which  
a {\it{diffusion}}, i.e.  having the continuous paths, corresponds. \\
{\textcolor{blue}{For the state space $S= {\cal S}'({\mathbb R}^d \to {\mathbb R})$, 
or subspaces of the tempered distributions ${\cal S}'({\mathbb R}^d \to {\mathbb R})$, 
the cases of {\it{infinite}} dimensional topological vector spaces, local Dirichlet forms have been formulated by, e.g., [AR 89, 90, 91], where the corresponding results are the constructions of 
${\cal S}'({\mathbb R}^d \to {\mathbb R})$-valued diffusion processes.
}}

{\textcolor{blue}{1-4. An example on Non-local (jump type Markov process ) 
Dirichlet form on real  $f(x)$, $x \in 
S \subset {\mathbb R}^m$, $m < \infty$, i.e. {\it{finite}} dimensional topological vector spaces
)}}
\medskip
\medskip

 Let us consider a toy model of $1$-dimensional {\it{non-local}} stochastic quantiztaon as follows: \, Let the state space $S \equiv \{- \frac12, \frac12\} \subset {\mathbb R}^1$.
Suppose that we are given some $p \in [0,1]$, and a probabilty measure $\mu$ on $S$ such that $ \mu(\{\frac12 \}) = p$, $\mu(\{- \frac12 \}) = 1-p$. Take $f, \,g \in {\cal S}({\mathbb R} \to {\mathbb R})$, and denote $u_1 = f(\frac12)$, $u_2 = f(- \frac12)$, $v_1 = g(\frac12)$ and $v_2 = g(- \frac12)$. Define an inner product $(f, g)_{L^2(\mu)}$ such that
$\displaystyle{(f, g)_{L^2(\mu)} \equiv \int_{\mathbb R} f(x) \, g(x) \, \mu(dx) 
= u_1 v_1 p + u_2 v_2 (1-p).}$\\

Then, $L^2(\mu)$ is the space of real valued bounded measurable functions on ${\mathbb R}$. On $L^2(\mu)$, define a closable  non-local Markovian symmetric form ${\cal E}_{(\alpha)}$, $\alpha \in (0,2)$ such that (see ($****$))
$$
{\cal E}_{(\alpha)} (f, g)  \equiv  \int_{x \ne y} \frac{(f(x)- f(y)) \, (g(x) - g(y))}{ |x- y|^{1 + \alpha}} \mu(dx) \, \mu (dy) $$
$$
= (u_1-u_2)(v_1-v_2) p (1-p) + (u_2 -u_1)(v_2 - v_1) (1-p)p 
$$
$$
= 
\left( A \left( \begin{array}{c}
                   u_1\\
                   u_2
                   \end{array}
             \right),
           \left( \begin{array}{c}
                   v_1\\
                   v_2
                   \end{array}
             \right)
\right)_{L^2(\mu)}
= \left( \left( \begin{array}{c}
                   u_1\\
                   u_2
                   \end{array}
             \right),
          A \left( \begin{array}{c}
                   v_1\\
                   v_2
                   \end{array}
             \right)
\right)_{L^2(\mu)},
$$
where 
$$ A = 2 \left( \begin{array}{cc}
    1-p & p-1 \\
    -p   & p 
                                    \end{array} \right)
.$$

The closed extension of ${\cal E}_{(\alpha)}$, denoted by the same notation, is a 
Dirichlet form with  domain ${\cal D}({\cal E}_{(\alpha)}) = L^2(\mu)$. 
 Then, $A$
 is the generator of Markovian semi group $\displaystyle{ e^{-A t}}$, $t \geq 0$ such that 
\begin{equation*} e^{-At} = \sum_{n= 0}^{\infty} \frac{(-1)^n t^n 2^n }{n !}  
\left( \begin{array}{cc}
    1-p & p-1 \\
    -p   & p 
                                    \end{array} \right)^n 
    \end{equation*}                              
$$= \left( \begin{array}{cc}
    p & 1-p \\
    p   & 1-p
                                    \end{array} \right)
+ e^{-2t}
\left( \begin{array}{cc}
    1-p & p-1 \\
    -p   & p 
                                    \end{array} \right).
$$
Denote the right hand side, the matrix, by $M_t$, $t \geq 0$, then it holds that
$$(p, 1-p) M_t = (p, 1-p).$$

Thus $M_t$, and hence $e^{-At}$ defines a continuous time conservative Markov process with  invariant measure  $\mu$, which is a non-local stochastic quantization of $\mu$. Obviously, for 
$L^2(\mu)$ it is impossible to formulate a corresponding local Dirichlet form. 
Roughly speaking, 
for a given random field $\mu$, 
in order to consider the stochastic quantizations through the arguments 
by means of the local Dirichlet forms, we have to suppose a sufficient regularity for $\mu$.

{\textcolor{blue}{There has been no general formulation on the Non-local Dirichlet forms defined on {\it{infinite}} dimensional toplogical vector spaces, before [AKYY2021] and [AKKYY2021]. 
The following is the rough explanation of the framework developed there.}}

\section{ Abstract theorems on non-local Dirichlet forms\\ on 
{\it{infinite}} dimensional topological vector spaces, \\
weighted $l^p$}
Let the state space $S$ be 
a weighted $l^p$ space, denoted by  $l^p_{(\beta_i)}$, 
such that, for some $p \in [1, \infty)$ and a weight $(\beta_i)_{i \in {\mathbb N}}$ with  $\beta_i \geq 0, 
i \in {\mathbb N}$,  
\begin{eqnarray}
S = l^p_{(\beta_i)} 
&\equiv& \big\{  {\mathbf x} = (x_1, x_2, \dots) \in {\mathbb R}^{\mathbb N} \, : 
\nonumber \\
&{}&
\| {\mathbf x} \|_{l^p_{(\beta_i)}} \equiv (  \sum_{i=1}^{\infty} {\beta}_i |x_i|^p )^{\frac1p} < \infty    \big\},
\end{eqnarray}
or 
a weighted $l^{\infty}$ space, denoted by  $l^{\infty}_{(\beta_i)}$, 
such that for  a weight $(\beta_i)_{i \in {\mathbb N}}$ with $\beta_i \geq 0, 
i \in {\mathbb N}$,  
\begin{eqnarray}
S = l^{\infty}_{(\beta_i)} 
&\equiv& \big\{  {\mathbf x} = (x_1, x_2, \dots) \in {\mathbb R}^{\mathbb N} \, : 
\nonumber \\
&{}&
\| {\mathbf x} \|_{l^{\infty}_{(\beta_i)}} \equiv \sup_{i \in {\mathbb N}} \beta_i |x_i| < \infty    \bigr\},
\end{eqnarray}

or 
\begin{equation}
S = {\mathbb R}^{\mathbb N},
\end{equation}
 the direct product space with the metric 
$d(\cdot, \cdot)$ such that for ${\mathbf x}, \,{\mathbf x}' \in {\mathbb R}^{\mathbb N}$, 
$d(\mathbf x, {\mathbf x}') \equiv \sum_{k=1}^{\infty} (\frac12)^k 
\frac{\|\mathbf x- {\mathbf x}'\|_k}{\|\mathbf x- {\mathbf x}'\|_k +1}$, 
 with 
$\| \mathbf x\|_k = ( \sum_{i=1}^k (x_i)^2)^{\frac12}, \, \,  {\mathbf x} =(x_1, x_2, \dots) \in {\mathbb R}^{\mathbb N}.
 $

Denote 
by ${\cal B}(S)$ the 
 Borel $\sigma$-field of $S$. 
Let $\mu$ be a given Borel probability measure on $(S, {\cal {B}}(S))$. 
For each $i \in {\mathbb N}$, 
{\textcolor{red}{let  $\sigma_{i^c}$  be the 
sub $\sigma$-field
}}
 of 
${\cal B}(S)$ that is generated by the Borel sets 
\begin{equation}
B = \left\{ {\mathbf x} \in S \, \, \, \big| \, x_{j_1} \in B_1, \dots x_{j_n} \in B_n \right\}, 
\end{equation}
$ j_k \ne i, \, \, B_k \in {\cal B}^1, \, \, k=1, \dots, n, \, \, n \in {\mathbb N}, \, \,
$
where ${\cal B}^1$ denotes the Borel $\sigma$-field of ${\mathbb R}^1$,

For each $i \in {\mathbb N}$, let $\mu(\cdot \, \big| \, \sigma_{i^c})$ be the conditional 
probability,
a one-dimensional probability distribution-valued  $\sigma_{i^c}$ measurable function, 
 that is characterized by
\begin{equation}
\mu \big( \{ {\mathbf x} \, \, : \, x_i \in A \} \cap B \big)
= \int_B \mu(A \, \big| \, {\sigma}_{i^c})   \, \mu(d {\mathbf x}),  
\end{equation}
$ \forall A \in {\cal B}^1, \, \,
\forall B \in {\sigma}_{i^c}.
$
Define
$$L^2(S; \mu) \equiv \Big\{ f \, \, \,  \Big| \, f : S \to {\mathbb R}, \, {\mbox{measurable and }}$$
\begin{equation}
\|f \|_{L^2} = \Bigl( \int_S |f({\mathbf x}) \mu(d {\mathbf x})  \Bigr)^{\frac12} < \infty \Big\},
\end{equation}
also define 
$$
{\cal F}C^{\infty}_0 \equiv 
{\mbox{the $\mu$ equivalence class of}} \,
$$
\begin{equation}
\left\{ f \, \, \, \Big| \, \exists n \in {\mathbb N}, \, f  \in C^{\infty}_0({\mathbb R}^n \to {\mathbb R}) \right\} \subset 
L^2(S; \mu).
\end{equation}

On $L^2(S;\mu)$, 
for $0 < \alpha <2$, define
the Markovian symmetric forms ${\cal E}_{(\alpha)}$  called  {\it{individually adapted Markovian symmetric form 
of index $\alpha$ to the measure $\mu$}}:\\
Firstly, 
for each $0 < \alpha <2$ and  $i \in {\mathbb N}$,
and
for 
the variables \,
$y_i, \, y'_i \in {\mathbb R}^1$, ${\mathbf x}  = (x_1, \dots , x_{i-1},x_i, x_{i+1}, \dots) \in S$ and  \,
${\mathbf x} \setminus x_i \equiv (x_1, \dots, x_{i-1},x_{i+1}, \dots)$, let 
\begin{eqnarray} 
\lefteqn{
\Phi_{\alpha}(u,v;y_i,y'_i,{\mathbf x} \setminus x_i)
} \nonumber \\
 &&\equiv 
\frac{1}{|y_i- y'_i|^{2\alpha +1}} 
\times 
\Big\{
u(x_1, \dots, x_{i-1}, y_i, x_{i +1}, \dots)
\nonumber \\
&& - u(x_1, \dots, x_{i-1}, y'_i, x_{i +1}, \dots) \Big\}
 \nonumber \\
&&\times 
\Big\{
v(x_1, \dots, x_{i-1}, y_i, x_{i +1}, \dots) 
\nonumber\\
&&- v(x_1, \dots, x_{i-1}, y'_i, x_{i +1}, \dots)
\Big\},
\end{eqnarray}


then, 
for each $ 0 < \alpha \leq 1$ and $i \in {\mathbb N}$, define
\begin{eqnarray}
\lefteqn{
{\cal E}_{(\alpha)}^{(i)}(u,v) 
}
\nonumber \\
&&\equiv 
\int_S
\Big\{
\int_{\mathbb R} I_{\{y_i \ne x_i\}}(y_i) \,
\Phi_{\alpha}(u,v;y_i,x_i, {\mathbf x}\setminus x_i)  
\nonumber \\
&&\times
 {\mu \big( dy_i \, \big| \, \sigma_{i^c} \big) }
\Big\} \mu(d{\mathbf x}).
\end{eqnarray}
and
\begin{equation}
{\cal E}_{(\alpha)}(u,v) \equiv \sum_{i \in {\mathbb N}} {\cal E}_{(\alpha)}^{(i)}(u,v),
\end{equation}
where 
$I_{\{\cdot\}}$ denotes the indicator function. 


For $y_i \ne y'_i$, 
(8) is well defined for any real valued ${\cal B}(S)$-measurable functions $u$ and $v$.
For {\it{the Lipschiz continuous}} functions 
${\tilde{u}} \in C^{\infty}_0({\mathbb R}^n \to {\mathbb R}) \subset {\cal F}C^{\infty}_0$ resp.
${\tilde{v}} \in C^{\infty}_0({\mathbb R}^m \to {\mathbb R}) \subset {\cal F}C^{\infty}_0$, \, $n,m \in {\mathbb N}$ which are representations of 
$u \in {\cal F}C^{\infty}_0$ resp.
$v \in {\cal F}C^{\infty}_0$, \, $n,m \in {\mathbb N}$, (9) and (10)  are well defined (the right hand side of (10) has only a finite number of sums).  In Theorem 1 given below we see that (9) and (10) are well defined for ${\cal F}C^{\infty}_0$, the space of $\mu$-equivalent class.

For  $1 < \alpha <2$,  we suppose that 
for each $i \in {\mathbb N}$, the conditional distribution 
 $\mu(\cdot \, \big| \, \sigma_{i^c})$ can be expressed by 
a locally bounded
 probability density $\rho(\cdot \, \big| \, \sigma_{i^c})$, 
$\mu$-a.e..
Precisely (cf. (2.5) of [AR91]), there exists a $\sigma_{i^c}$-measurable function 
{\textcolor{blue}{
$0 \leq \rho(\cdot \, \big| \, \sigma_{i^c})$ on ${\mathbb R}^1$ 
}}
and
\begin{equation}
\mu(dy \, \big| \, \sigma_{i^c}) = 
\rho(y \, \big| \, \sigma_{i^c}) \, dy, \quad \mu-a.e.,
\end{equation}
holds, with  $\rho(\cdot \, \big| \, \sigma_{i^c})$ a function such that for any compact 
$K \subset {\mathbb R}$ there exists an $L_i <\infty$, which may depend on $i$,  and for any $y \in K$,
{\textcolor{blue}{
\begin{equation}
{\rm{ess}} \, \sup_{y \in {\mathbb R}^1}  \rho(y \, \big| \, \sigma_{i^c}) \leq L_i, \qquad  \mu-a.e.,
\end{equation}
where ${\rm{ess}}\sup_{y \in {\mathbb R}^1}$ is taken with respect to the Lebesgue measure on ${\mathbb R}^1$.
}}
Then define 
the non-local form ${\cal E}_{(\alpha)}(u,v)$, for  $1 < \alpha <2$, by the same formula as (10).

\begin{remark}
{\textcolor{blue}{
{\it{
For the ${\cal B}(S)$ measurable function 
$\int_{\mathbb R} 1_{\{y_i \ne x_i \}} \Phi_{\alpha}(u,v;y_i,x_i, {\mathbf x} \setminus x_i)
\mu(dy_i \, | \, \sigma_{i^c})$ 
by taking the expectation conditioned by the sub $\sigma$-field $\sigma_{i^c}$, 
it holds that 
(cf.,[Fukushima, Uemura 2012]):
\begin{eqnarray}
\lefteqn{
{\cal E}_{(\alpha)}^{(i)}(u,v) } \nonumber \\
&&\equiv
\int_S
\Big\{
\int_{{\mathbb R}} I_{\{y_i \ne x_i\}}\,
\Phi_{\alpha}(u,v;y_i,x_i, {\mathbf x \setminus x_i})  \, {\mu \big( dy_i \, \big| \, \sigma_{i^c} \big)} 
\Big\} \mu(d{\mathbf x}) \nonumber \\
&&=
\int_S
\int_{\mathbb R} \Big\{
\int_{{\mathbb R}} I_{\{y_i \ne x_i\}}\,
\Phi_{\alpha}(u,v;y_i,x_i, {\mathbf x \setminus x_i})   
\nonumber \\
&&\times
 {\mu \big( dy_i \, \big| \, \sigma_{i^c} \big)}
\Big\}
 {\mu \big( dx_i\, \big| \, \sigma_{i^c} \big) }
\, \mu(d{\mathbf x}) \nonumber \\
&&=
\int_S
\Big\{
\int_{{\mathbb R}^2} I_{\{y_i \ne y'_i\}}\,
\Phi_{\alpha}(u,v;y_i,y'_i, {\mathbf x \setminus x_i})  
\nonumber \\
&&\times
{\mu \big( dy_i \, \big| \, \sigma_{i^c} \big)}\,
 {\mu \big( dy'_i \, \big| \, \sigma_{i^c} \big) }
\Big\} \mu(d{\mathbf x})   \nonumber
\end{eqnarray}
}}
}}
\end{remark}


\begin{theorem}
{\textcolor{blue}{
{\it{
The symmetric non-local forms ${\cal E}_{(\alpha)}$, $0 <\alpha <2$ given by (10) 
(for $1 < \alpha <2$ with the additional assumption (11) 
with (12)
) are 
(cf. Remark 1-i),ii))
\\
i) \quad \, \, well-defined on ${\cal F}C^{\infty}_0$; \\
ii) \quad \, Markovian;\\
iii) \quad closable in $L^2(S;\mu)$. \\
}}
}}
\end{theorem}

\begin{remark}
{\textcolor{blue}{
{\it{
 For 
$1 < \alpha < 2$ the assumption (11) with (12) can be replaced by the following general one:\,  for each compact $K \subset {\mathbb R}$, there exists an 
$L_i < \infty$ and 
\begin{equation*}
\sup_{y \in {\mathbb R}} I_K(y) \int_{{\mathbb R}} \frac{I_K(y')}{|y -y'|^{2 \alpha -1}} \mu(dy'\,|\, 
\sigma_{i^c})  \, \leq L_i, \quad \mu-a.e..
\end{equation*}
}}
}}
\end{remark}

\begin{remark}
{\textcolor{red}{
A way to exclude the trivial cases.}} 
Suppose that $\mu$ is not a trivial probability measure, i.e., there exists no ${\mathbf x} \in {\mathbb R}^{\mathbb N}$ such that $\mu({\mathbf x}) =1$. Eventhough it might happen that ${\cal E}_{\alpha} \equiv 0$, and the resulting Markov prosess becoms trivial, i.e., it does not have any probabilistic dynanics. By the following (i) or (ii), such cases can be excluded:\\
{\textcolor{blue}{
(i) \, \, 
{\it{
For $0 < \alpha \leq 1$, 
in order to exclude the case where ${\cal E}_{\alpha} \equiv  0$, or 
equivaently the case where the associated Markov process $X_t = {\mathbf x}$, $\forall t \geq 0$, for some ${\mathbf x} \in S$, 
we should put the assumption such that there exists an $ {\mathbf x } \in S$. an $i \in {\mathbb N}$, an 
 $L >0$ and 
$$\mu(\{{\mathbf y} \, : \, y_i \in (x_i-L, x_i + L)\} \, \big| \, {\mathbf x} \setminus x_i) >0.$$
Note  that the  (11) and (12) are sufficient to this.
}}
}}

{\textcolor{blue}{
(ii) \, \, Under tha assuption that $\mu$ is not trivial, there must exist $i, \, j \in {\mathbb N}$, two of distinguished points $(x_i^0, x^0_j), \, ({x'}_i^0, {x'}_j^0) \in {\mathbb R}^2$, for e.g., $x^0_i \ne {x'}^0_i$,  and $\epsilon >0$ such that 
$$\mu( |x_i - x^0_i| < \epsilon \, \big| \, x_i = \frac{x_j^0 - {x'}_j^0}{x_i^0 - {x'}_i^0} 
+ b) > 0,$$
$$\mu( |x_i - {x'}^0_i| < \epsilon \, \big| \, x_i = \frac{x_j^0 - {x'}_j^0}{x_i^0 - {x'}_i^0} 
+ b) > 0,$$
where $ b = x^0_j -  \frac{x_j^0 - {x'}_j^0}{x_i^0 - {x'}_i^0} \cdot x_i^0$. 
For such $i$, in the definition (9) and (10), if we replace $\mu(dy_i|{\sigma}_{i^c})$ by 
$\mu(dy_i | y_i = \frac{x_j^0 - {x'}_j^0}{x_i^0 - {x'}_i^0} 
+ b, {\sigma}_{i^c})$, then the symmetric form with this modification becomes non-trivial one.}}
{\textcolor{red}{
(END of Remark 3.)
}}
\end{remark}

%
%


To prove the Theorem 1, we have to show that \\
i-1) \quad 
for $u \in {\cal F}C^{\infty}_0$ such that $u=0, \, \mu-a.e.$, it holds that\quad 
$\displaystyle{{\cal E}_{(\alpha)}(u,u) =0}$, and \\
i-2) \quad for any $u, v \in {\cal F}C^{\infty}_0$, \, ${\cal E}_{(\alpha)}(u,v) \in {\mathbb R}$, \\
For the statement ii), we have to show that (cf. [Fukushima]) for any $\epsilon >0$ there exists a real function $\varphi_{\epsilon}(t)$, $ -\infty < t < \infty$, such that ${\varphi}_{\epsilon}(t) = t, \, \forall t \in [0,1]$, 
$-\epsilon \leq {\varphi}_{\epsilon}(t) \leq 1 + \epsilon, \, \forall t \in (-\infty, \infty)$, and
$0 \leq {\varphi}_{\epsilon}(t')-{\varphi}_{\epsilon}(t) \leq t'-t$ for $t < t'$, such that for any $u \in {\cal F}C^{\infty}_0$ it holds that ${\varphi}_{\epsilon}(u) \in {\cal F}C^{\infty}_0$ and
\begin{equation}
{\cal E}_{(\alpha)}({\varphi}_{\epsilon}(u), {\varphi}_{\epsilon}(u)) \leq {\cal E}_{(\alpha)}(u,u).
\end{equation}

For the statement iii), we have to show the following: \, For a sequence $\{u_n\}_{n \in {\mathbb N}}$, 
$u_n \in {\cal F}C^{\infty}_0$, $n \in {\mathbb N}$, if 
\begin{equation}
\lim_{n \to \infty} \| u_n \|_{L^2(S; \mu)} = 0,
\end{equation}
and
\begin{equation}
\lim_{n,m \to \infty} {\cal E}_{(\alpha)}(u_n -u_m, u_n -u_m) =0,
\end{equation}
then
\begin{equation}
\lim_{n \to \infty} {\cal E}_{(\alpha)} (u_n,u_n) =0.
\end{equation}

Proof of i-1): \, 
For each $i \in {\mathbb N}$ and any
real valued ${\cal B}(S)$-measurable function $u$, note that for each $\epsilon >0$,
$$  I_{\{\epsilon < |x_i-y_i|\} }(y_i) \, I_K(y_i)
\Phi_{\alpha}(u,u;y_i, x_i, {\mathbf x} \setminus x_i)$$
defines a ${\cal B}(S \times {\mathbb R})$-measurable function.
The function \\
$\Phi_{\alpha}(u,u;y_i,x_i, {\mathbf x} \setminus x_i)$, is defined. 
 ${\cal B}(S \times {\mathbb R})$ is the Borel $\sigma$-field of $S \times {\mathbb R}$.
${\mathbf x} = (x_i, i \in {\mathbb N}) \in S$ and $y_i \in {\mathbb R}$.
Then, 
for any compact subset $K$ of ${\mathbb R}$, 
$0 \leq I_{\{\epsilon < |x_i-y_i|\} }(y_i) \, I_K(y_i)
\Phi_{\alpha}(u,u;y_i, x_i, {\mathbf x} \setminus x_i)$ converges monotonically to
$ I_{\{y_i \ne x_i\}}(y_i) \,
\Phi_{\alpha}(u,u;y_i,x_i, {\mathbf x} \setminus x_i)$ as $ K \uparrow {\mathbb R}$ and $\epsilon \downarrow 0$, for every $y_i \in {\mathbb R}$, \, ${\mathbf x} \in S$, 
and by the
{\textcolor{blue}{
Fatou's Lemma}}, 
we have

\begin{eqnarray}
\lefteqn{
\int_S  
\left\{
\int_{\mathbb R} I_{\{y_i \ne x_i\}}(y_i) \,
\Phi_{\alpha}(u,u;y_i,x_i, {\mathbf x} \setminus x_i)  \, {\mu \big( dy_i \, \big| \, \sigma_{i^c} \big) }
\right\} 
} \nonumber \\
&&\times \mu(d{\mathbf x})
  \nonumber \\
&& = 
\int_S  \liminf_{K \uparrow {\mathbb R}} \liminf_{\epsilon \downarrow 0} 
 \Big\{ 
 \int_{\mathbb R}  I_{\{\epsilon < |x_i-y_i|\} }(y_i) \, I_K(y_i) 
\nonumber \\
&&\times
\Phi_{\alpha}(u,u;y_i, x_i, {\mathbf x} \setminus x_i)  \, {\mu \big( dy_i \, \big| \, \sigma_{i^c} \big) }
\Big\} \mu(d{\mathbf x}) 
\nonumber \\
&& \leq 
\liminf_{K \uparrow {\mathbb R}} 
\liminf_{\epsilon \downarrow 0}
\int_S   
 \Big\{ 
 \int_{\mathbb R}  I_{\{\epsilon < |x_i-y_i| \}}(y_i) \, I_K(y_i) 
\nonumber \\
&&\times
\Phi_{\alpha}(u,u;y_i, x_i, {\mathbf x} \setminus x_i)  \, {\mu \big( dy_i \, \big| \, \sigma_{i^c} \big) }
\Big\} \mu(d{\mathbf x}), 
\end{eqnarray}

where $K$ denotes a compact set of  ${\mathbb R}$. 
For any $\epsilon >0$,  
\begin{eqnarray}
\lefteqn{
\int_S \Big\{ \int_{\mathbb R}  I_{\{\epsilon < |x_i-y_i| \}}(y_i) \, I_K(y_i) \, 
\frac{1}{|y_i -x_i|^{2 \alpha + 1}} 
}
\nonumber \\
&&\times 
\big(u(x_1, \dots) \big)^2 
\, {\mu \big( dy_i \, \big| \, \sigma_{i^c} \big) }
\Big\} \mu(d{\mathbf x})
 \nonumber \\
&&\leq 
\frac{1}{{\epsilon}^{2 \alpha +1}} 
\int_S \Big\{ \int_{\mathbb R}  I_{\{\epsilon < |x_i-y_i| \}}(y_i) \, I_K(y_i) 
\nonumber \\
&&\times 
 \big(u(x_1, \dots) \big)^2 
\, {\mu \big( dy_i \, \big| \, \sigma_{i^c} \big) }
\Big\} \mu(d{\mathbf x}) 
 \nonumber \\
&&\leq 
\frac{1}{{\epsilon}^{2 \alpha +1}} 
\int_S \Big\{ \int_{\mathbb R}  
 \big(u(x_1, \dots) \big)^2 
\nonumber \\
&&\times 
\, {\mu \big( dy_i \, \big| \,  \sigma_{i^c} \big) }
\Big\} \mu(d{\mathbf x}) 
  \\
&&=
\frac{1}{{\epsilon}^{2 \alpha +1}} 
\int_S 
 \big(u(x_1, \dots) \big)^2 
 \mu(d{\mathbf x})  = 0, \nonumber
\end{eqnarray}

also,
$$
\int_S 
\big(u(x_1, \dots) \big)^2 
\Big\{ \int_{\mathbb R}  I_{\{\epsilon < |x_i-y_i| \}}(y_i) \, I_K(y_i) 
$$
$$ \times 
\frac{1}{|y_i -x_i|^{2 \alpha + 1}} 
\, {\mu \big( dy_i \, \big| \, \sigma_{i^c} \big) }
\Big\} \mu(d{\mathbf x})
$$
\begin{equation}
\leq
\frac{1}{{\epsilon}^{2 \alpha +1}} 
\int_S 
 \big(u(x_1,\dots) \big)^2 
 \mu(d{\mathbf x})  = 0,
\end{equation}
and 
from (18), by the Cauchy Schwaz's inequality
\begin{eqnarray}
\lefteqn{ 
\Big|
\int_S 
u(x_1, \dots) 
\Big\{ \int_{\mathbb R}  I_{\{\epsilon < |x_i-y_i| \}}(y_i) \, I_K(y_i) 
}
\nonumber \\
&&\times \frac{1}{|y_i -x_i|^{2 \alpha + 1}} 
u(x_1, \dots)  
\, {\mu \big( dy_i \, \big| \, \sigma_{i^c} \big) }
\Big\} \mu(d{\mathbf x}) \Big| 
\nonumber
\\
&&\leq
\frac{1}{{\epsilon}^{2 \alpha +1}} 
\int_S 
 \big(u(x_1, \dots) \big)^2 
 \mu(d{\mathbf x}),
\end{eqnarray}

By (8) and (17), from (18), (19) and (20),
for any Borel measurable function $u$ on $S$ such that 
$$
u(x_1, \dots) = 0, \qquad \mu-a.e.,
$$ 
it holds that
{\textcolor{blue}{
$${\cal E}_{(\alpha)}^{(i)} (u,u) = 0, \quad   \forall i \in {\mathbb N},
\qquad {\cal E}_{(\alpha)} (u,u) = 0.$$
}}

{\textcolor{blue}{
In order to show i-2), 
 {
{for $0 < \alpha \leq 1$}},  
}}
take any 
{\it{representation}} 
${\tilde{u}} \in C^{\infty}_0({\mathbb R}^n)$ of $u \in {\cal F}C^{\infty}_0$, $n \in {\mathbb N}$. 
 {
{Using  $ 0 < \alpha + 1 \leq 2$, it is easy to see  from the definition (8) that}}
there exists an $M < \infty$ depending on ${\tilde{u}}$  such that
\begin{equation}
0 \leq \Phi_{\alpha}({\tilde{u}},{\tilde{u}};y_i,y'_i,{\mathbf x} \setminus x_i) \leq M,
 \quad {\mbox{$ \forall {\mathbf x} \in S$,  and $\forall y_i, \, y'_i \in {\mathbb R}$}}.
\end{equation}
Since, $u = {\tilde u} + {\overline{0}}$ for some 
real valued ${\cal B}(S)$-measurable function ${\overline{0}}$ such that ${\overline{0}}=0$, $\mu$-a.e., 
by (21) together with i-1) and the the Cauchy Schwarz's inequality, for $u \in {\cal F}C^{\infty}_0$, ${\cal E}_{(\alpha)}(u,u) \in {\mathbb R}$, $0 < 
\alpha \leq 1$,  is identical with 
${\cal E}_{(\alpha)}({\tilde{u}},{\tilde{u}})$ and 
well-defined 
 (in fact, for only a finite number of $i \in {\mathbb N}$. 
we have 
${\cal E}^{(i)}_{(\alpha)}(u,u)  \ne 0$, cf. also (10)). 
Then by the Cauchy Schwarz's inequality i-2) follows.

A proof of ii) is very similar to the one given in section 1 of [Fukushima].

Proof of 
iii): \,  \, Suppose that a sequence $\{u_n \}_{n \in {\mathbb N}}$ satisfies (14) and (15).  Then, by (14) there exists a measurable set ${\cal N} \in {\cal B}(S)$ and a 
{
{
sub sequence $\{ u_{n_k} \}$ of $\{u_n\}$ such that
$$
\mu({\cal N}) = 0, \quad   \lim_{{n_k} \to \infty} u_{n_k} (\mathbf x) = 0, \quad \forall {\mathbf x} \in S \setminus {\cal N}.
$$
Define
$$
{\tilde{u}}_{n_k} ({\mathbf x}) = u_{n_k}({\mathbf x}) \, {\mbox{for ${\mathbf x} \in S \setminus {\cal N}$, \quad and}} \quad 
{\tilde{u}}_{n_k} ({\mathbf x}) = 0 \, {\mbox{for ${\mathbf x} \in  {\cal N}$}}.
$$
Then, 
\begin{equation} {\tilde{u}}_{n_k} (\mathbf x) = u_{n_k}(\mathbf x),  \, \, \mu-a.e., 
\qquad  \lim_{n_k \to \infty} {\tilde{u}}_{n_k} (\mathbf x) = 0, \quad \forall {\mathbf x} \in S.
\end{equation}
}}

By the fact  i-1),  for each $i$,
$$
\int_S  
\Big\{
 \int_{\mathbb R} I_{\{y_i \ne x_i\}}(y_i) \,
\Phi_{\alpha}(u_n,u_n;y_i,x_i, {\mathbf x} \setminus x_i)  \, {\mu \big( dy_i \, \big| \, \sigma_{i^c} \big) }
\Big\} \mu(d{\mathbf x})
$$
$$
= \int_S  
 \Big\{ 
\int_{\mathbb R} I_{\{y_i \ne x_i\}}(y_i) \,
\lim_{n_k \to \infty}
\Phi_{\alpha}(u_n - {\tilde{u}}_{n_k},u_n- {\tilde{u}}_{n_k};y_i,x_i, {\mathbf x} \setminus x_i)  
$$
$$ \times  {\mu \big( dy_i \, \big| \, \sigma_{i^c} \big) }
\Big\} \mu(d{\mathbf x})
$$
$$\leq
\liminf_{n_k \to \infty}
\int_S 
\Big\{ 
\int_{\mathbb R} I_{\{y_i \ne x_i\}} \,
\Phi_{\alpha}(u_n - {\tilde{u}}_{n_k},u_n- {\tilde{u}}_{n_k};y_i, x_i, {\mathbf x} \setminus x_i)  
$$
$$\times
 {\mu \big( dy_i \, \big| \, \sigma_{i^c} \big) }
\Big\} \mu(d{\mathbf x})
$$
$$=
\liminf_{n_k \to \infty}
\int_S  \, 
 \Big\{ 
 \int_{\mathbb R} I_{\{y_i \ne x_i\}}\,
\Phi_{\alpha}(u_n - {{u}}_{n_k},u_n- {{u}}_{n_k};y_i, x_i, {\mathbf x} \setminus x_i)  $$
$$\times 
{\mu \big( dy_i \, \big| \, \sigma_{i^c} \big) }
\Big\} \mu(d{\mathbf x})
$$
\begin{equation}
\equiv
\liminf_{n_k \to \infty}
{\cal E}^{(i)}_{(\alpha)}(u_n -u_{n_k}, u_n -u_{n_k}).
\end{equation}

Now, by using the assumption (15) to the right hand side of (23), we get 
\begin{equation}
\lim_{n \to \infty} {\cal E}_{(\alpha)}^{(i)}(u_n,u_n) = 0, \qquad \forall i \in {\mathbb N}.
\end{equation}
{\textcolor{blue}{
(24) together with i) shows that for each $i \in {\mathbb N}$, ${\cal E}^{(i)}_{(\alpha)}$ with the domain ${\cal F}C_0^{\infty}$ is closable in $L^2(S;\mu)$.
Since, ${\cal E}_{(\alpha)} \equiv \sum_{i \in {\mathbb N}} {\cal E}^{(i)}_{(\alpha)}$, by using the Fatou's Lemma, from (24) and the assumption (15) we see that 
$$
{\cal E}_{(\alpha)}(u_n,u_n) = \sum_{i \in {\mathbb N}} \lim_{m \to \infty} {\cal E}^{(i)}_{(\alpha)}(u_n -u_m, u_n -u_m) $$
$$
\leq \liminf_{m \to \infty}  {\cal E}_{(\alpha)}(u_n -u_m, u_n -u_m)  \to 0 \, {\mbox{as $n \to 
\infty$}}.
$$
This proves  (16) 
(cf.  Proposition I-3.7 of [MR] for a general argument). The proof of iii) is completed.
}}

{\textcolor{blue}{
{\bf{The proof of Theorem 1,  for $1< \alpha < 2$}}
}}

The proof of i-1), ii) and  iii)   can be carried out by the 
 completely 
same manner as the previous proof 
we have provided 
for the case $0 < \alpha \leq 1$. 
We  only  show that  i-2), i.e., ${\cal E}_{(\alpha)}(u,u) < \infty, \, \forall u \in {\cal F}C^{\infty}_0$ 
also holds when we make use of 
 the additional assumption (11) with (12),.\\

The detailed proof is omitted.

For each $i \in {\mathbb N}$, denote by $X_i$ the random variable that represents the coordinate $x_i$ of  ${\mathbf x} = (x_1, x_2, \dots)$:
\begin{equation}
X_i \, : S \ni {\mathbf x} \longmapsto x_i \in {\mathbb R}.
\end{equation}
Then
\begin{equation}
\int_S 1_{B}(x_i) \, \mu(d {\mathbf x}) = \mu (X_i \in B), \quad {\mbox{for \, \, $B \in 
{\cal B}(S)$}}.
\end{equation}

\begin{theorem}
{(Strictly Quasi-regularlity)} 
{\textcolor{blue}{
{\it{
Let \,  $0 < \alpha \leq 1$, and 
$({\cal E}_{(\alpha)}, {\cal D}({\cal E}_{(\alpha)}))$ be the closed Markovian symmetric form defined  through Theorem 1.\\
i) \quad In case where $S= l^p_{(\beta_i)}$, $ 1 \leq p < \infty$, if there exists
positive $l^p$ sequence $\{\gamma^{-1}_i \}_{i \in {\mathbb N}}$ ( for e.g., $\gamma_i = i^{\frac{1 + \delta}p}$ for some $\delta >0$), and 
 an \, $M_0 < \infty$ \, and 
{\textcolor{blue}{
\begin{equation}
\sum_{i=1}^{\infty} 
\beta^{\frac{2}p}_i \gamma_i^2 \cdot
\mu \Big( |X_i| > M_0 \cdot \beta_i^{-\frac{1}{p}}\, \gamma_i^{-1}\Big) < \infty,
\end{equation}
\begin{equation}
\mu \Big( \bigcup_{M \in {\mathbb N}} \{ |X_i| \leq M \cdot \beta_i^{-\frac{1}{p}}\, \gamma_i^{-\frac{1}p}, \, \forall i \in  {\mathbb N} \}  \Big) =1,
\end{equation}
}}
hold, then $({\cal E}_{(\alpha)}, {\cal D}({\cal E}_{(\alpha)}))$ is a strictly quasi-regular Dirichlet form. 
}}
}}

{\textcolor{blue}{
{\it{
ii) \quad 
 In case where $S= l^{\infty}_{(\beta_i)}$ defined by (2), if there exist an \, $M_0 < \infty$ and a  
sequence $\{ \gamma_i \}_{i \in {\mathbb N}}$ such that $0 < \gamma_1  \leq \gamma_2 \leq \cdots \to \infty$, and  both 
{\textcolor{blue}{
\begin{equation}
 \sum_{i=1}^{\infty} 
\beta_i^2 {\gamma}_i^2 \cdot
\mu \Big( |X_i| > M_0 \cdot \beta_i^{-1}\, \gamma_i^{-1} \Big) < \infty,
\end{equation}
\begin{equation}
\mu \Big( \bigcup_{M \in {\mathbb N}} \{  |X_i| \leq M \cdot \beta_i^{-1}\, \gamma_i^{-1}, \, \forall i \in  {\mathbb N} \}  \Big) =1,
\end{equation}
}}
hold, then $({\cal E}_{(\alpha)}, {\cal D}({\cal E}_{(\alpha)}))$ is a strictly quasi-regular Dirichlet form. \\
iii) \quad 
 In case where $S= {\mathbb R}^{\mathbb N}$ defined by (3), 
 $({\cal E}_{(\alpha)}, {\cal D}({\cal E}_{(\alpha)}))$ is a strictly quasi-regular Dirichlet form.
}}
}}
\end{theorem}

Proof of Theorem 2. \quad We have to show that 
$({\cal E}_{(\alpha)}, {\cal D}({\cal E}_{(\alpha)}))$ satisfies \\
i) \quad There exists an ${\cal E}_{(\alpha)}$-nest $(D_M)_{M \in {\mathbb N}}$ consisting of compact sets. \\
ii) \quad There exists a subset of ${\cal D}({\cal E}_{(\alpha)})$, that is dense with respect to the norm $\displaystyle{\| \cdot \|_{L^2(S;\mu)} + \sqrt{{\cal E}_{(\alpha)}}}$. 
And the elements of the subset have ${\cal E}_{(\alpha)}$-quasi continuous versions. \\
iii) \quad There exists $u_n \in {\cal D}({\cal E}_{(\alpha)})$, $n \in {\mathbb N}$, having ${\cal E}_{(\alpha)}$-quasi continuous $\mu$-versions ${\tilde{u}}_n$, $n \in {\mathbb N}$, and an ${\cal E}_{(\alpha)}$-exceptional set ${\cal N} \subset S$ such that $\{{\tilde{u}}_n \, :\, 
n \in {\mathbb N} \}$ separates the points of $S \setminus {\cal N}$.\\
iv) \quad For the {\it{strictly}} quasi-regularity, it suffices to show that 
$1 \in {\cal D}({\cal E}_{(\alpha)})$.

For the case where 
$S = l^p_{(\beta_i)}$, 
for simplicity, let $ \gamma_i^{-1} = i^{-\frac{1 + \delta}p}$ for some 
$\delta >0$. A key point of the proof is the fact that  for each $M \in {\mathbb N}$, 
\begin{equation}
D_M \equiv \left\{ {\mathbf x} \in l^p_{(\beta_i)} \, : \, \, {\beta}_i^{\frac{1}{p}} | x_i|  \leq 
M \cdot i^{-\frac{1 + \delta}p}, \,  i \in {\mathbb N} 
 \right\},
\end{equation}
is a compact set in $S = l^p_{(\beta_i)}$. \\
{\textcolor{blue}{Note that $D_M$ is not 
identical to but a proper subset of the bounded set such that 
$$ \Big\{ {\mathbf x}  \in  l^p_{(\beta_i)} \, : \,  \big( \sum \beta_i |x_i|^p \Big)^{\frac1p} \leq M \big\}.$$
}}

Let \, $\eta(\cdot) \in C^{\infty}_0({\mathbb R})$ be a function such that \, $\eta(x) \geq 0$, \, 
$\displaystyle{|\frac{d}{dx} \eta(x)| \leq 1}$, \, $\forall x \in {\mathbb R}$ \, and 
\begin{equation}
\eta(x) = \left\{ 
           \begin{array}{ll}
             1, & \qquad \qquad  |x| \leq 1; \\ [0.2cm]
             0, & \qquad \qquad  |x| \geq 3.
           \end{array} \right.
\end{equation}
For each $M \in {\mathbb N}$\, and \,$i \in {\mathbb N}$, let 
$
\eta_{M,i}(x) \equiv \eta \left(M^{-1} \cdot i^{\frac{1+ \delta}p} \beta_i^{\frac{1}p}\cdot x \right), \quad x \in {\mathbb R},
$\\
 then,
 $\displaystyle{\prod_{i \geq 1} \eta_{M,i} \in l^p_{(\beta_i)}}$, \,
$
{\rm{supp}} \, \big[ \prod_{i \geq 1}\eta_{M,i} \big] \subset D_{3M}$, \,  $M \in {\mathbb N}$.\\
For each $f \in C^{\infty}_0({\mathbb R}^n \to {\mathbb R})$, \, $n \in {\mathbb N}$, define
\begin{equation}
f_M(x_1, \dots, x_n, x_{n+1}, \dots)  \equiv f(x_1, \dots, x_n) \cdot \prod_{i \geq 1} \eta_{M,i}(x_i). \end{equation}

Under the condition (27), 
it is possible to show that $f_M \in {\cal D}({\cal E}_{(\alpha)})$. 
Also, by (28), it is possible to show that there exists a subsequence $\{f_{M_l} \}_{l \in {\mathbb N}}$ of $\{f_M\}_{M \in {\mathbb N}}$ such that the {\it{Cesaro mean}}
$$
w_m \equiv \frac{1}m \sum_{l=1}^m f_{M_l} \to u\equiv f \cdot \prod_{i \geq 1} 1_{{\mathbb R}}(x_i)$$ 
\begin{equation}
{\mbox{in \quad ${\cal D}({\cal E}_{(\alpha)})$ \quad as \quad $n\to \infty$}}.
\end{equation}
(34) shows that 
$$
{\mbox{the linear hull of}} \,
\Big\{ f_M, \, M \in {\mathbb N}\, \, : \, 
f \in C^{\infty}_0({\mathbb R}^n \to {\mathbb R}), \, n \in {\mathbb N}  \Big\}.
$$
can be taken as an ${\cal D}({\cal E}_{(\alpha)})$-nest.

\begin{theorem}
{(Strictly Quasi-regularlity)}
{\textcolor{blue}{
{\it{
Let \,  
$1< \alpha < 2$. 
 Suppose that the assumption (11) with (12) hold.
Let  $({\cal E}_{(\alpha)}, {\cal D}({\cal E}_{(\alpha)}))$ be the closed  Markovian symmetric form defined at the beginning of this section through Theorem 1.
Then the following statements hold:\\
i)  \quad In the case where $S= l^p_{(\beta_i)}$, $ 1 \leq p < \infty$, as defined by (1),  if there exists a
positive $l^p$ sequence $\{\gamma^{-\frac{1}p}_i \}_{i \in {\mathbb N}}$  and 
 an  $M_0 < \infty$,  and both (28),
\begin{equation}
\sum_{i = 1}^{\infty} \big(\beta^{\frac{1}p}_i {\gamma}_i^{\frac{1}p} \big)^{ \alpha +1} \cdot
\mu \Big(\beta_i^{\frac{1}p} |X_i| > M_0 \cdot {\gamma}_i^{-\frac{1}{p}} \Big) < \infty,
\end{equation}
\begin{equation}
\lim_{M \to \infty} 
M^{-{\alpha}} 
\sum_{i=1}^{\infty} L_{M,i} \cdot
\big(\beta^{\frac1p}_i {\gamma}_i^{\frac{1}p} \big)^{ \alpha}\cdot
\mu \Big( \beta_i^{\frac{1}P} |X_i| > M \cdot {\gamma}_i^{-\frac{1}{p}} \Big) < \infty,
\end{equation}
}}
}}

{\textcolor{blue}{
{\it{
hold, then $({\cal E}_{(\alpha)}, {\cal D}({\cal E}_{(\alpha)}))$ is a strictly quasi-regular Dirichlet form,
where for each  
$M \in {\mathbb N}$ and $i \in {\mathbb N}$,
 $L_{M,i}$ is the bound of the conditional probability density $\rho$ for a given compact set  
$$K_{M,i} \equiv \big[ -6M\cdot \beta_i^{-\frac1p} \, {\gamma}_i^{-\frac{1}p}, 6M\cdot
\beta_i^{-\frac1p} \,  {\gamma}_i^{-\frac{1}p} \big] \subset {\mathbb R}$$
in the assumption (12).
}}
}}

{\textcolor{blue}{
{\it{
ii) \quad 
 In the case where $S= l^{\infty}_{(\beta_i)}$ as defined by (2), if there exists a 
sequence $\{ \gamma_i \}_{i \in {\mathbb N}}$ such that $0 < \gamma_1  \leq \gamma_2 \leq \cdots \to \infty$, and  both 
\begin{equation}
 \sum_{i=1}^{\infty} 
\big(\beta_i {\gamma}_i \big)^{ \alpha +1} \cdot
\mu \Big( \beta_i |X_i| > M_0 \cdot \gamma_i^{-1} \Big) < \infty, \quad {\mbox{for some $M_0 < \infty$}},
\end{equation}
\begin{equation}
\lim_{M \to \infty} M^{-{\alpha}} \sum_{i=1}^{\infty} L_{M,i} \cdot
\big(\beta_i {\gamma}_i \big)^{ \alpha} \cdot
\mu \Big( \beta_i |X_i| > M \cdot  \gamma_i^{-1} \Big) < \infty,
\end{equation}
and (30) hold, then $({\cal E}_{(\alpha)}, {\cal D}({\cal E}_{(\alpha)}))$ is a strictly quasi-regular Dirichlet form,
 where for each  $M \in {\mathbb N}$ and $i \in {\mathbb N}$,
 $L_{M,i}$ is the bound of the conditional probability density $\rho$ for a given compact set  
$$K_{M,i} \equiv \big[ -6M\cdot \beta_i^{-1} \, \gamma_i^{-1}, 6M\cdot
\beta_i^{-1} \,  \gamma_i^{-1} \big] \subset {\mathbb R}$$
in the assumption (12).
}}
}}

{\textcolor{blue}{
{\it{
iii) \quad 
 In the case where $S= {\mathbb R}^{\mathbb N}$,  if there exists a sequence 
$\{ \gamma_i \}_{i \in {\mathbb N}}$ such that $0 < \gamma_i$, $\forall i \in {\mathbb N}$, and 
\begin{equation}
\lim_{M \to \infty} M^{-{\alpha}} \sum_{i = 1}^{\infty} L_{M,i} \cdot \gamma_i^{-{\alpha}} \cdot 
\mu \Big( |X_i| > M \cdot \gamma_i \Big) < \infty,
\end{equation}
holds, then 
 $({\cal E}_{(\alpha)}, {\cal D}({\cal E}_{(\alpha)}))$ is a strictly quasi-regular Dirichlet form, 
where for each  $M \in {\mathbb N}$ and $i \in {\mathbb N}$,
 $L_{M,i}$ is the bound of the conditional probability density $\rho$ for a given compact set  
$$K_{M,i} \equiv \big[ -6M\cdot \gamma_i, \, 6M\cdot
\gamma_i \big] \subset {\mathbb R}.$$
}}
}}
\end{theorem}

\begin{theorem}
{\textcolor{blue}{
{\it{
Let $0  < \alpha <2$, and let 
${\cal E}_{(\alpha)}, {\cal D}({\cal E}_{(\alpha)}))$ be a strictly quasi-regular Dirichlet form on 
$L^2(S; \mu)$ that is defined through Theorem 2 or Theorem 3.
Then for $({\cal E}_{(\alpha)}, {\cal D}({\cal E}_{(\alpha)}))$, there exists a properly associated $\mu$-tight special standard process, in short  a strong Markov process taking values in $S$ and having right continuous trajectories with left limits up to the life time (cf. Definitions IV-1.5, 1.8 and 1.13 of [MR] for its
precise definition), 
$${\mathbb M} \equiv \Big(\Omega, {\cal F}, (X_t)_{t \geq 0}, (P_{\mathbf x})_{\mathbf x \in 
S_{\triangle}} \Big),$$
where $\triangle$ is an adjoined extra points, called as the cemetery, of $S$.\
}}
}}
\end{theorem}

\section{Applications and Examples on infinite dimensional situations}
Euclidean  ({\it{scalar}}) quantum fields are expressed as random fields on ${\cal S}' \equiv {\cal S}'({\mathbb R}^d \to {\mathbb R})$, or, resp., 
${\cal S}'({\mathbb T}^d \to {\mathbb R})$, the Schwartz's space of real tempered distributions on ${\mathbb R}^d$, resp.,  the $d$-dimensional torus ${\mathbb T}^d$, with $d \geq 1$ a given space time dimension.  Hence, each Euclidean quantum field is taken as a probability space $\displaystyle{ \big( {\cal S}', 
{\cal B}({\cal S}'), \nu \big)}$, where ${\cal B}({\cal S}')$ is the Borel $\sigma$-field of ${\cal S}'$ and $\nu$ is a Borel probability measure on ${\cal S}'$. One of the standard 
theorem through which such $\nu$ are constructed is the Bochner-Minlos's theorem 
(see below), which is an existence theorem of probability measures on 
Hilbert nuclear spaces.  Since the space ${\cal S}$ and its dual ${\cal S}'$ is a Hilbert nuclear space, and by making use of a Hilbert-Schmidt operators defined  on it, we can adapt our Theorems.

\subsection{ How to apply the abstract theorems to concrete probability\\ measures
 on function spaces (spaces of paths)}
As a most ideal model, consider the case where the state space is the Sobolev space 
such that 
$H^1_0 = H^1_0( [0, 2 \pi]   \to {\mathbb R})$ namely, 
the closure of $\{ f \, : \, f \in C^1([0, 2 \pi] \to {\mathbb R}), \, f(0) = f(2 \pi) = 0 \}$ with respect to the norm, by the {\it{Pincare inequality}}, 
$\| f \| = \|\frac{d}{dx}f \|_{L^2([0,2 \pi] \to {\mathbb R})}$. \\
Let $\{ a_n \}_{n \in {\mathbb N}}$, be the Fourier coefficient of $f \in H^1_0$, i.e., 
$$a_n = \frac{1}{\pi} \int_{0}^{2 \pi} f(x) \cdot \sin nx \, dx, \quad n = 1,2, \dots.$$

Then the map $\tau$, 
{\textcolor{blue}{ the
identification between $H^1_0$ and  
$l^2_{\{n^2 \}}$ 
}}
such that 
$$\tau \, : H^1_0 \ni f \longmapsto (a_1, a_2, \dots ) \in l^2_{\{n^2 \}},$$
defines an {\textcolor{blue}{
{\it{isometric isomorphism}}
}}
, where $l^2_{\{n \}}$ is the weighted $l^2$ space with the norm 
$\| {\mathbf x} \|_{l^2_{\{n \}} }
= ( \sum_{i=1}^{\infty} n^2 |a_i|^2 )^{\frac12}$, for 
${\mathbf x} = (a_1, a_2, \dots) \in l^2_{\{n \}}$.  Recall that by the Fourier expansion, 
{\textcolor{blue}{
$$f(x) = \sum_{n=1}^{\infty} a_n \sin nx.$$
}}
Suppose that we are given a probability measure $\nu$ on $H^1_0$, then define 
$$\mu (B) = \nu (\{f \, : \, \tau(f) \in B \}), \quad {\mbox{for $B \in {\cal B}(l^2_{\{n \}})$ }}.
$$

Then, by applying the abstract theorems to the state space $S = l^2_{\{n^2 \}}$ with the probability measure $\mu$, if we have an objective Markov process $(X)_{t \geq 0}$ taking the values in $l^2_{\{n^2 \}}$ of which invariant measure is $\mu$, we then define an $H^1_0$-valued Markov process $(Y_t)_{t \geq 0}$ such that
{\textcolor{blue}{
$$Y_t = \sum_{n=1}^{\infty} x_n(t) \, \sin nx, \quad {\mbox{for}} \, \,  X_t = ((x_1(t), x_2(t), \dots)).$$
}}
The function valued, i.e., the path space valued, Markov process  
$(Y_t)_{t \geq 0}$ is the stochastic quantization of $\nu$, which we want in the concrete physical situation. 
{\textcolor{blue}{
Namely, $X_t$ is the Fourier coefficients valued Markov process, and $Y_t$ is the corresponding function space, $H^1_0$, valued Markov process.
}}

\subsection{3-2. \, Bochner-Minlos's theorem}

Let us recall the Bochner-Minlos theorem stated in a general framework.  Let $E$ be a nuclear space ( cf., e.g., Chapters 47-51 of [Tr{\`e}ves 67]). Suppose in particular that $E$ is a countably Hilbert space,   characterized by a sequence of {\it{real}} Hilbert norms $\|{}\, \,{}\|_n$, $n \in {\mathbb N}\cup \{0\}$ such that 
$\|{}\, \,{}\|_0 < \|{}\, \,{}\|_1 < \cdots <\|{}\, \,{}\|_n < \cdots$.
Let $E_n$ be the completion of $E$ with respect to the norm $\|{}\, \,{}\|_n$, then by definition 
$E = \bigcap_{n \geq 0} E_n$ and $E_0 \supset E_1 \supset \cdots \supset E_n \supset \cdots$. Define 
$$E^{\ast}_n \equiv {\mbox{the dual space of $E_n$, and assume 
the identification}}$$  $$E^{\ast}_0 = E_0.$$

Then  we have 
$$
E \subset \cdots \subset E_{n+1} \subset E_n \subset \cdots \subset  E_0$$
$$ = E^{\ast}_0 \subset 
\cdots \subset E^{\ast}_n \subset E^{\ast}_{n+1} \subset \cdots \subset E^{\ast}.
$$
Since by assumption $E$ is a nuclear space, for any $m \in {\mathbb N}\cup\{0\}$ there exists an $n 
\in {\mathbb N}\cup \{0\}$, $n > m$, such that the (canonical) injection $T^n_m\, : E_n \to E_m$ is a trace class (nuclear class) positive operator. The Bochner-Minlos theorem 
 (cf. [Hida 80])
is given as follows:

\begin{theorem}{\bf{(Bochner-Minlos Theorem)}} \\
Let $C(\varphi)$, $\varphi \in E$, be a complex valued function on $E$ such that\\
i) \quad \,\, $C(\varphi)$ is continuous with respect to the norm $\| \, \cdot \, \|_{m}$ for some $m \in {\mathbb N} \cup \{0\}$; \\
ii) \quad \, ({\bf{positive definiteness}}) \quad 
for any $k \in {\mathbb N}$,
$$
\sum_{i,j =1}^k {\bar{\alpha}}_i \alpha_j C({\varphi}_i - {\varphi}_j) \geq 0, \,  
\forall \alpha_i \in {\mathbb C}, \, \, \forall {\varphi}_i \in E, \,\,  i=1,\dots, k;
$$
(where ${\bar{\alpha}}$ means complex conjugate of $\alpha$).

\it{
iii) \quad ({\bf{normalization}}) \quad  $C(0) = 1$. \\
Then, there exists a unique Borel probability measure $\nu$ on $E^{\ast}$ such that 
$$
C(\varphi) = \int_{E^{\ast}} e^{i <\phi, \varphi>} \nu(d \phi), \qquad \varphi \in E.
$$
Moreover,
if the (canonical) injection $T^n_m\, : E_n \to E_m$, for all $ n > m$, is a Hilbert-Schmidt operator, then 
the support of $\nu$ is in $E^{\ast}_n$, where 
$<\phi, \varphi> = {}_{E^{\ast}}<\phi, \varphi>_{E}$ is the dualization between $\phi \in E^{\ast}$ and 
$\varphi \in E$.
}
\end{theorem}

\subsection{3-3. \, Applications of Hilbert Schmidt Theorem (Isomorphism of weighted $l^2$)}
 Let
\begin{equation*}
{\cal H}_0 \equiv \Big\{ f \, : \, \|f\|_{{\cal H}_0} = \big((f,f)_{{\cal H}_0} \big)^{\frac12} < \infty, \, \, f : {\mathbb R}^d \to {\mathbb R}, 
\end{equation*}
\begin{equation}
 {\mbox{measurable}} \Big\} \supset 
{\cal S}({\mathbb R}^d),
\end{equation}
where
\begin{equation}
(f,g)_{{\cal H}_0} \equiv (f,g)_{L^2({\mathbb R}^d)} = \int _{{\mathbb R}^d} f(x) g(x) \, dx.
\end{equation}
Let
\begin{equation}
H \equiv (|x|^2 + 1)^{\frac{d+1}2} (- \Delta +1)^{\frac{d +1}2} (|x|^2 + 1 )^{\frac{d +1}2},
\end{equation}
\begin{equation}
H^{-1} \equiv (|x|^2 + 1)^{- \frac{d+1}2} (- \Delta +1)^{- \frac{d +1}2} (|x|^2 + 1 )^{- 
\frac{d +1}2},
\end{equation}
be the pseudo differential operators on

${\cal S}'({\mathbb R}^d \to {\mathbb R}) \equiv {\cal S}'({\mathbb R}^d)$
with the $d$-dimensional Laplace operator $\Delta$. 
For each $n \in {\mathbb N}$, define 
\begin{equation*}
{\cal H}_n \equiv {\mbox{the completion of ${\cal S}({\mathbb R}^d)$
with respect to 
the norm }}
\end{equation*}
\begin{equation}
\|f \|_n  = \sqrt{(f,f)_n} \, \,  \mbox{ with} \, \,  
(f,g)_n = (H^n f, H^n g) _{{\cal H}_0},
\end{equation}
and
\begin{equation*}
{\cal H}_{-n} \equiv {\mbox{the completion of ${\cal S}'({\mathbb R}^d)$ with respect to 
the norm }}
\end{equation*}
\begin{equation}
\|f \|_{-n}
= \sqrt{(f,f)_{-n}} \, \, \mbox{with} \, \,  
(f,g)_{-n} = ((H^{-1})^n f, (H^{-1})^n g) _{{\cal H}_0}.
\end{equation}
by taking an inductive limit ${\cal H} = \bigcap_{n \in {\mathbb N}} {\cal H}_n$, 
then
\begin{equation*}
{\cal H} \subset \cdots \subset {\cal H}_{n + 1} \subset {\cal H}_n \subset \cdots \subset  {\cal H}_0 \subset  \cdots \subset {\cal H}_{-n} \subset {\cal H}_{-n-1}
\end{equation*}
\begin{equation}
 \subset \cdots \subset {\cal H}^{\ast}.
\end{equation}

For the positive self-adjoint operator $H^{-1}$ on ${\cal H}_0 = L^2({\mathbb R}^d \to {\mathbb R})$,take the orthonormal base (O.N.B.) $\{{\varphi}_i\}_{i \in {\mathbb N}}$ of ${\cal H}_0$ such that 
\begin{equation}
H^{-1} {\varphi}_i = \lambda_i \, {\varphi}_i, \qquad i \in {\mathbb N},
\end{equation}
where $\{\lambda_i\}_{i \in {\mathbb N}}$ is the corresponding eigenvalues such that 
$1 \geq {\lambda}_1 \geq \lambda_2 \geq \cdots > 0$, which satisfies 
\begin{equation}
\sum_{i \in {\mathbb N}} (\lambda_i)^2 < \infty, \quad {\mbox{i.e.,}} \quad 
\{\lambda_i\}_{i \in {\mathbb N}} \in l^2.
\end{equation}
Then,
\begin{equation}
\{(\lambda_i)^{n} {\varphi}_i \}_{i \in {\mathbb N}} \quad {\mbox{is an O.N.B. of ${\cal H}_n$ }}
\end{equation}
and 
\begin{equation}
\{(\lambda_i)^{-n} {\varphi}_i \}_{i \in {\mathbb N}} \quad {\mbox{is an O.N.B. of ${\cal H}_{-n}$ }}
\end{equation}

Thus, 
by the Fourier series expansion 
for $f \in {\cal H}_m$,
$$
f  = \sum_{i \in {\mathbb N}} a_i (\lambda_i^{m} {\varphi}_i), \quad {\mbox{with}} 
$$
\begin{equation}
a_i \equiv \big(f, \, (\lambda_i^{m} {\varphi}_i) \big)_{m} = {\lambda_i^{-m}} (f, \,
\varphi_i)_{L^2}, 
\, \, i \in {\mathbb N},
\end{equation}
we have an 
{\textcolor{blue}{
{\it{isometric isomorphism}} 
 $\tau_m$
 }}
  for each $m \in {\mathbb Z}$ such that 
\begin{equation}
\tau_m \, : \, {\cal H}_m \ni f \longmapsto ({\lambda}_1^m a_1, {\lambda}_2^m a_2, \dots) \in l^2_{(\lambda_i^{-2m})},
\end{equation}
where  $l^2_{(\lambda_i^{-2m})}$ is the weighted $l^2$ space defined by (1) with $p=2$, 
and  $\beta_i  = \lambda_i^{-2m}$.

\begin{equation}
{\cal H}_2 \, \, \subset \, \,  \, \, {\cal H}_1 \, \, \subset {\cal H}_0 = L^2({\mathbb R}^d) \subset {\cal H}_{-1} \subset {\cal H}_{-2},
\end{equation}

\begin{equation}
 l^2_{(\lambda_i^{-4})} \, \, \subset  l^2_{(\lambda_i^{-2})} \, \, \subset \, \, \, \, \, \, \, \, l^2 \, \, \, \, \, \, \, \, \,
\subset  \, \, \, \, \, \, \,  l^2_{(\lambda_i^2)} \subset  \, \, \, \, \, \, \,  l^2_{(\lambda_i^{4})}.
\end{equation}
(53) and (54) show the correspondences between the function spaces and the weighted $l^2$ spaces through ${\tau}_m$.

{\bf{Example 0. (The Euclidean free quantum field)}} \quad

Let ${\nu}_0$ be the Euclidean free field 
probability 
measure on ${\cal S}' \equiv {\cal S}'({\mathbb R}^d)$.
It is characterized by 
the (generalized) characteristic function  $C(\varphi)$  in Theorem 4 of $\nu_0$  given by  
\begin{equation*}
C(\varphi) = \exp ({- \frac12 (\varphi, (- \Delta + m^2_0)^{-1} \varphi)_{L^2({\mathbb R}^d)}}).
\end{equation*}
$$
 {\mbox{for \, $\varphi \in {\cal S}({\mathbb R}^d \to {\mathbb R})$}},
$$

Equivalently, $\nu_0$ is the centered Gaussian probability measure on ${\cal S}'$, the covariance of which 
  is given by
\begin{equation*}
\int_{{\cal S}'} <\phi, {\varphi}_1>   <\phi, {\varphi}_1> \, \nu_0(d \phi) 
= \big(\varphi_1, (- \Delta + m^2_0)^{-1} \varphi_2 \big)_{L^2({\mathbb R}^d)}, 
\end{equation*}
$$
 \varphi_1, \, \varphi_2 \in {\cal S}({\mathbb R}^d \to {\mathbb R}),
$$
where $\Delta$ is the $d$-dimensional Laplace operator and $m_0 >0$

( for $d \geq 3$, we can 
 also 
allow for   $m_0 =0$) is a given mass for this  scalar field.
$\phi (f) = <\phi, f>$, $f \in {\cal S}({\mathbb R}^d \to {\mathbb R})$ is 
 the coordinate process $\phi$ to $\nu_0$ (for the Euclidean free field cf. 
[Pitt 71], 
[Nelson 74], 
 and, e.g., 
[Simon 74], [Glimm,Jaffe 87], 
[A,Y 2002], [A,Ferrario,Y 2004]). The functional $C(\varphi)$ is continuous with respect to the  norm of the space  ${\cal H}_0 =L^2({\mathbb R}^d)$, and the kernel of 
$(- \Delta + m^2_0)^{-1}$, which is the Fourier inverse transform of $(|\xi|^2 + m^2_0)^{-1}$, $\xi \in {\mathbb R}^d$, is explicitly given by Bessel functions (cf., e.g., section 2-5 of  [Mizohata 73]).
Then, by 
Bochner-Minlos's theorem
the support of $\nu_0$ 
can be taken to be 
 in the separable Hilbert  spaces  ${\cal H}_{-n}$, $ n \geq 1$.

Let us apply Theorems 1, 2 and 3 with $p = \frac12$ to this random field.

 We start the  consideration from the case where $\alpha= 1$, a simplest situation. 
 Then, we shall state the corresponding results for the cases where $0 < \alpha <1$.
 Now, we take ${\nu}_0$ as a Borel probability measure on ${\cal H}_{-2}$. By  taking $m = -2$, 
$\tau_{-2}$ defines an 
{\textcolor{blue}{
isometric isomorphism
}}
 such that 
\begin{equation*}
\tau_{-2} \, : \, {\cal H}_{-2} \ni f \longmapsto (a_1, a_2, \dots) \in l^2_{(\lambda_i^{4})},
\end{equation*}
$$
 {\mbox{with}} \quad 
a_i \equiv (f, \, \lambda_i^{-2} {\varphi}_i)_{-2}, \, \, \, i \in {\mathbb N}.
$$
Define a probability measure $\mu$ on 
$ l^2_{(\lambda_i^{4})}$ such that 
$$\mu(B) \equiv \nu_0 \circ \tau^{-1}_{-2}(B) \quad {\mbox{for}} \quad B \in {\cal B}(l^2_{(\lambda_i^{4})}).
$$
We set $S = l^2_{(\lambda_i^{4})}$ in Theorems 1, 2 and 3, then it follows that the weight $\beta_i$ satisfies $\beta_i = \lambda_i^4$.

 We can take ${\gamma_i}^{-\frac12} = \lambda_i$ in Theorem 2,  then, we have
\begin{equation}
\sum_{i=1}^{\infty} 
\beta_i \gamma_i \cdot
\mu \Big( \beta_i^{\frac12} |X_i| > M \cdot \gamma_i^{-\frac{1}2}\Big) 
\leq \sum_{i =1}^{\infty} \beta_i \gamma_i = \sum _{i =1}^{\infty} (\lambda_i)^2 < \infty
\end{equation}
(55) shows that the condition (27) holds.

Also, as has been mentioned above, 
since $\nu_0({\cal H}_{-n}) = 1$,  for any   $ n \geq 1$, we have 
$$1 = \nu_0({\cal H}_{-1}) = \mu(l^2_{(\lambda_i^2)}) 
$$
$$= 
\mu \big( \bigcup_{M \in {\mathbb N}} \{ |X_i| \leq M \beta_i^{- \frac12} \gamma_i^{- \frac12}, \, \forall i \in {\mathbb N} \} \big),
$$
$$
 {\mbox{for $\beta_i = \lambda_i^4$, \, $\gamma_i^{- \frac12} = \lambda_i$}}.
$$
This shows that the condition (28) is satisfied.

Thus, by Theorem 2 and Theorem 3,  
for   $\alpha = 1$,
there exists an  $l^2_{((\lambda_i)^{4})}$-valued  Hunt process 
\begin{equation}
{\mathbb M} \equiv \big(\Omega, {\cal F}, (X_t)_{t \geq 0}, (P_{\mathbf x})_{\mathbf x \in 
S_{\triangle}} \big),
\end{equation}
associated  to the non-local Dirichlet form
$({\cal E}_{(\alpha)}, {\cal D}({\cal E}_{(\alpha)}))$.

We can now define an
${\cal H}_{-2}$-valued 
 process 
$(Y_t)_{t \geq 0}$ 
such that
\begin{equation} 
(Y_t)_{t \geq 0} \equiv \big({\tau}^{-1}_{-2}(X_t) \big)_{t \geq 0}.
\end{equation}
Equivalently,  for $X_t = (X_1(t), X_2(t), \dots) \in l^2_{(\lambda_i^4)}$, $P_{\mathbf x}-a.e.$, by setting $A_i(t)$ such that 
 $A_i(t) \equiv \lambda_i X_i(t)$,  we see that  
$Y_t$ \, is also given by 
\begin{equation*}
Y_t = \sum_{i \in {\mathbb N}} A_i(t) (\lambda_i^{-2} \varphi_i) = \sum_{i \in {\mathbb N}} X_i(t) \varphi_i \in {\cal H}_{-2},  \forall t \geq 0, \, \, P_{\mathbf x}-a.e., 
\end{equation*}
$$
{\mbox{for any $x \in S_{\triangle}$}}.
$$

Then, $Y_t$ is an ${\cal H}_{-2}$-valued Hunt process that is a {\it{stochastic quantization}} (according to the definition we gave to this term) with respect to the 
non-local Dirichlet form 
$({\tilde{\cal E}}_{(\alpha)}, {\cal D}({\tilde{\cal E}}_{(\alpha)}))$ 
on $L^2({\cal H}_{-2}, \nu_0)$, that is defined through
$({\cal E}_{(\alpha)}, {\cal D}({\cal E}_{(\alpha)}))$, by making use of $\tau_{-2}$. 
This holds for  $\alpha = 1$.   
\medskip
\medskip

{\textcolor{blue}{
For the cases where $0 < \alpha <1$, we can also apply Theorems 1, 2 and 3, and then have the corresponding result to (57).}}

{\bf{Example 1. (The $\Phi^4_3$ Euclidean free quantum field)}} \quad 
An alternative procedure, through which the $\Phi^4_d$ Euclidean field measures are defined, is the following: \, Let $d= 2,3$. For each bounded region $\Lambda \subset 
{\mathbb R}^d$ and a lattice spacing $\epsilon >0$, let 
\begin{equation}
L_{\epsilon, \Lambda} \equiv (\epsilon {\mathbb Z})^d \cap \Lambda,
\end{equation}
 and define a family of real valued random variables $\phi \equiv \{ \phi(x) \, : \, x \in L_{\epsilon, \Lambda} \}$, the probability distribution of which is given by
\begin{equation}
\nu_{\epsilon, \Lambda} \equiv \frac{1}{Z_{\epsilon, \Lambda}} \prod_{x \in L_{\epsilon, \Lambda}} e^{- S_{L_{\epsilon, \Lambda}} (\phi)} d \phi(x),
\end{equation}
where $Z_{\epsilon, \Lambda}$ is the normalizing constant, and

\begin{equation*}
S_{L_{\epsilon, \Lambda}} (\phi) \equiv \frac12 \sum_{<x,y>} \epsilon^{d-2} \big(\phi(x) -\phi(y) \big)^2 
\end{equation*}
$$
+ \frac12 a_{\epsilon} \sum_{x \in L_{\epsilon, \Lambda}} \epsilon^d \phi^2(x)
+ \frac{\lambda}2 \sum_{x \in L_{\epsilon, \Lambda}} \epsilon^d \phi^4(x),
$$
with $a_{\epsilon}$ a counter term depending on $\epsilon >0$ and $d=2, 3$, $\lambda \geq 0$ a coupling constant;
  $<x,y>$ denotes the nearest neighbor points in $L_{\epsilon, \Lambda}$. 
In [Brydges,Fr{\"o}hlich,Sokal 83] it is shown, roughly speaking
that, for adequately small $\lambda \geq 0$,  there exists a subsequence 
$\{ \nu_{\epsilon_i, \Lambda_j}\}_{i,j \in {\mathbb N}}$ 
of $\{ \nu_{\epsilon, \Lambda}\}_{\epsilon >0, \Lambda \subset {\mathbb R}^d}$ 
with $\lim_{i \to \infty} \epsilon_i = 0$ and $\lim_{j \to \infty} \Lambda_j = {\mathbb R}^d$, 
and a weak limit 
\begin{equation}
\nu^{\ast} \equiv \lim_{i \to \infty} \lim_{j \to \infty} \nu_{\epsilon_i, \Lambda_j}
\end{equation}
exists in the space of Borel probability measures

 on ${\cal S}'({\mathbb R}^d \to {\mathbb R})$ by interpreting $\nu_{\epsilon, \Lambda}$ as 
an element in this space,
 for each $\epsilon >0$ and $\Lambda \subset 
{\mathbb R}^d$  
For each $\epsilon >0$, and bounded region $\Lambda \subset {\mathbb R}^d$, $d=1,2,3$, and for $F(\phi)$ a polynomial in $\{ \phi(x) \, : \, x \in L_{\epsilon, \Lambda} \}$,  let 
\begin{equation}
< F >_{\epsilon, \Lambda} \equiv \int_{{\mathbb R}^{N(\epsilon, \Lambda)}} F(\phi) \nu_{\epsilon, \Lambda}(d \phi), 
\end{equation}
where 
$N(\epsilon, \Lambda)$ is the cardinality of $L_{\epsilon,\Lambda}$.
By [Sokal 82] (cf. also section 2 of [Brydges,Fr{\"o}hlich,Sokal 83]), the following limit exists:
\begin{equation}
<F>^{(\epsilon)} \equiv \lim_{\Lambda \uparrow {\mathbb R}^d} < F >_{\epsilon, \Lambda}.
\end{equation}
Also, for each $\epsilon >0$, and $d=1,2,3$, there exists a weak limit 
\begin{equation}
\nu_{\epsilon} \equiv \lim_{\Lambda \uparrow {\mathbb R}^d} \nu_{\epsilon, \Lambda},
\end{equation}

that is a Borel probability measure on ${\cal S}'({\mathbb R}^d \to {\mathbb R})$.  
Define 
\begin{equation}
S^{(\epsilon)} (x-y) \equiv <\phi(x) \cdot \phi(y)>^{(\epsilon)}, \qquad x, y \in (\epsilon {\mathbb Z})^d,
\end{equation}
and 
\begin{equation*}
S^{(\epsilon)}_n(x_1, \dots, x_n) \equiv < \prod_{i=1}^n \phi(x_i) >^{(\epsilon)}, \quad x_i \in L_{\epsilon}, 
\end{equation*}
$$
i=1, \dots, n, \, \, n \in {\mathbb N}.
$$

From Theorem 5 we see that 
the 
supports of the  Borel probability measure $\nu_{\epsilon}$ on ${\cal S}'({\mathbb R}^d \to {\mathbb R})$, for each $\epsilon >0$,  and of a weak limit of 
a subsequence of $\{\nu_{\epsilon} \}_{\epsilon >0}$, denoted by $\nu$, are 
 all 
in the Hilbert space 
${\cal H}_{-2} \subset {\cal S}'({\mathbb R}^d \to {\mathbb R})$ defined by (43), (44) and (45) 
so that from now on 
\begin{equation}
{\mbox{$\nu_{\epsilon}$ and $\nu$ can be understood as probability measures on ${\cal H}_{-2}$}}.
\end{equation}

For each $n \in {\mathbb N}$, we can take a sequence $\{\epsilon_{n,i} \}_{i \in {\mathbb N}}$, with 
$\epsilon_{n,i} >0$ and $\lim_{i \to \infty} \epsilon_{n,i} =0$ such that
$\{ {\cal S}^{(\epsilon_{n,i})}_{2n} \}_{i \in {\mathbb N}}$ converges weakly (as $i \to \infty$) to some $S_{2n} \in ({\cal H}_{-1})^{\otimes 2n}$.  By taking subsequences and using a diagonal argument, we then see that there exists a sequence $\{ \epsilon_i \}_{i \in {\mathbb N}}$ with $\epsilon_i > 0$ and 
$\lim_{i \to \infty} \epsilon_i = 0$ such that $\{ S^{\epsilon_i}_{2n} \}_{i \in {\mathbb N}}$ converges weakly (as $i \to \infty$) to $S_{2n} \in ({\cal H}_{-1})^{\otimes 2n}$  for any $n \in {\mathbb N}$.

By these, we can define
 the (generalized) characteristic function 
\begin{equation}
C(\varphi) \equiv \sum_{n =0}^{\infty} \frac{(-1)^n}{(2n)!}  \, \big< S_{2n}, \, {\varphi}^{\otimes 2n} \big>,
\end{equation}
that satisfies
\begin{equation}
|C(\varphi)| \leq e^{\frac12 K \|\varphi \|_{{\cal H}_1}}, \qquad \forall \varphi \in {\cal H}_1,
\end{equation}
and, for any $n \geq 2$
\begin{equation}
{\mbox{$\nu$ is the  probability measure on ${\cal H}_{-n}$ corresponding to $C(\varphi)$}}.
\end{equation}

By taking $n = -3$, 
$\tau_{-3}$ defines an 
{\textcolor{blue}{
isometric isomorphism 
}}
such that 
$$
\tau_{-3} \, : \, {\cal H}_{-3} \ni f \longmapsto (a_1, a_2, \dots) \in l^2_{(\lambda_i^{6})},
\quad {\mbox{with}} 
$$
\begin{equation} 
a_i \equiv (f, \, \lambda_i^{-3} {\varphi}_i)_{-3}, \, \, \, i \in {\mathbb N}.
\end{equation}
Define  $\mu$ on 
$ l^2_{(\lambda_i^{6})}$ such that 
\begin{equation}
\mu(B) \equiv \nu \circ \tau^{-1}_{-3}(B) \quad {\mbox{for}} \quad B \in {\cal B}(l^2_{(\lambda_i^{6})}).
\end{equation}
Set $S = l^2_{(\lambda_i^{6})}$. We can take $\beta_i = \lambda_i^6$,  ${\gamma_i}^{-\frac12} = \lambda_i$ in Theorem 2-i) 
with $p =2$,  then, 
\begin{equation}
\sum_{i=1}^{\infty} 
\beta_i \gamma_i \cdot
\mu \Big( \beta_i^{\frac12} |X_i| > M \cdot \gamma_i^{-\frac{1}2}\Big) 
\leq \sum_{i =1}^{\infty} \beta_i \gamma_i = \sum _{i =1}^{\infty} (\lambda_i)^4 < \infty.
\end{equation}

This shows that (27) holds, 
also, it is possible to see that (28) holds.\\
Thus, by Theorem 2-i) and Theorem 4,  
for each  $0 < \alpha \leq 1$,
there exists an  $l^2_{(\lambda_i^{6})}$-valued  Hunt process 
\begin{equation}
{\mathbb M} \equiv \big(\Omega, {\cal F}, (X_t)_{t \geq 0}, (P_{\mathbf x})_{\mathbf x \in 
S_{\triangle}} \big), 
\end{equation}
associated  to the non-local Dirichlet form
$({\cal E}_{(\alpha)}, {\cal D}({\cal E}_{(\alpha)}))$.
Then define an
${\cal H}_{-3}$-valued 
 process 
$(Y_t)_{t \geq 0}$ 
such that 
$(Y_t)_{t \geq 0} \equiv \big({\tau}^{-1}_{-2}(X_t) \big)_{t \geq 0}.$

Equivalently, for $X_t = (X_1(t), X_2(t), \dots) \in l^2_{(\lambda_i^6)}$, $P_{\mathbf x}-a.e.$, by setting $A_i(t)$ such that 
$X_i(t) = \lambda_i^{-3} A_i(t)$, then 
$Y_t$ \, is given by 
$$
Y_t = \sum_{i \in {\mathbb N}} A_i(t) (\lambda_i^{-3} \varphi_i) = \sum_{i \in {\mathbb N}} X_i(t) \varphi_i \in {\cal H}_{-3},$$
\begin{equation}
 \forall t \geq 0, \, \, P_{\mathbf x}-a.e..
\end{equation}
It is an ${\cal H}_{-3}$-valued Hunt process that is a {\it{stochastic quantization}} with respect to the 
non-local Dirichlet form 
$({\tilde{\cal E}}_{(\alpha)}, {\cal D}({\tilde{\cal E}}_{(\alpha)}))$ 
on $L^2({\cal H}_{-3}, \nu)$, that is defined through
$({\cal E}_{(\alpha)}, {\cal D}({\cal E}_{(\alpha)}))$, by making use of $\tau_{-3}$.





\end{document}